\documentclass{amsart}
\usepackage{amsmath}
\usepackage{amssymb}
\newtheorem{thm}{Theorem}[section]
\newtheorem{cor}[thm]{Corollary}
\newtheorem{conj}[thm]{Conjecture}
\newtheorem{lem}[thm]{Lemma}
\newtheorem{prop}[thm]{Proposition}
\newtheorem{cons}[thm]{Construction}
\theoremstyle{definition}
\newtheorem{defn}[thm]{Definition}
\theoremstyle{remark}
\newtheorem{rem}[thm]{Remark}
\newtheorem{ex}[thm]{Example}
\newtheorem{exs}[thm]{Examples}

\long\def\Cor#1{\begin{cor} #1 \end{cor}}

\long\def\Lem#1{\begin{lem} #1 \end{lem}}
\long\def\Prop#1{\begin{prop} #1 \end{prop}}

\long\def\Rem#1{\begin{rem} #1 \end{rem}}
\long\def\Ex#1{\begin{ex} #1 \end{ex}}

\def\bar#1{\overline{#1}}
\def\Sect{\section}
\def\Rarr#1#2{\xrightarrow[#2]{#1}}

\long\def\Ref#1#2#3#4#5#6{
\bibitem{#1}
{\rm #2,}
\textit{#3.}
{\rm #4}
\textbf{#5}
{\rm #6.}
}
\long\def\Refb#1#2#3#4{
\bibitem{#1}
{\rm #2,}
\textit{#3.}
#4.
}
\def\Zz{{\mathbb Z}}
\def\Rr{{\mathbb R}}
\def\Cc{{\mathbb C}}
\def\Ff{{\mathbb F}}
\def\Qq{{\mathbb Q}}
\def\Tt{{\mathbb T}}
\def\i{{\rm i}}
\def\phi{\varphi}
\def\into{\hookrightarrow}
\def\leq{\leqslant}
\def\geq{\geqslant}
\def\st{\mid}
\def\map{{\rm map}}
\def\Zero{{\rm Zero}}
\def\o{{\rm o}}
\def\aa{\mathfrak{a}}
\def\O{{\rm O}}
\def\Hom{{\rm Hom}}
\def\e{{\rm e}}
\def\U{{\rm U}}

\def\SS{\mathfrak{S}}
\begin{document}

\title{Borsuk--Ulam theorems for elementary abelian $2$-groups}

\author{M.~C.~Crabb}
\address{%
Institute of Mathematics\\
University of Aberdeen \\
Aberdeen AB24 3UE \\
UK}
\email{m.crabb@abdn.ac.uk}
\date{January 2022}
\begin{abstract}
Let $G$ be a compact Lie group and let $U$ and $V$ be 
finite-dimensional real $G$-modules
with $V^G=0$. A theorem of Marzantowicz,
de Mattos and dos Santos estimates the covering dimension of the
zero-set of a $G$-map from the unit sphere in $U$ to $V$ when
$G$ is an elementary elementary abelian $p$-group for some
prime $p$ or a torus.
In this note, the classical Borsuk--Ulam theorem will be used to give
a refinement of their result estimating the dimension
of that part of the zero-set on which an elementary abelian
$p$-group $G$ acts freely or a torus $G$ acts with finite isotropy
groups. 
The methods also provide an easy answer
to a question raised in \cite{DM}.  
\end{abstract}
\subjclass[2010]{55M20, 
55M25, 
55R25, 
55M35, 
55N91} 
\keywords{Borsuk--Ulam theorem, equivariant mapping, 
Euler class}
\maketitle
\par\noindent
Let $G$ be a compact Lie group 
and let $U$ and $V$ be finite-dimensional real $G$-modules
(which we may assume to be equipped with a $G$-invariant 
Euclidean inner product) with the fixed subspace $V^G$ equal to
zero. 
A theorem \cite[Theorem 2.1]{MdMdS} of Marzantowicz,
de Mattos and dos Santos estimates the covering dimension of the
zero-set of a $G$-map from $S(U)$, the unit sphere in $U$, to $V$
when $G$ is an elementary abelian $p$-group for some prime $p$
or a torus.
In this note, the classical Borsuk--Ulam theorem will be used to give
a refinement of their result which estimates the dimension
of that part of the zero-set on which an
elementary abelian $p$-group $G$ acts freely
or a torus $G$ acts with only finite isotropy groups. 

Section 1 reviews the relevant Borsuk--Ulam theorems. 
Elementary abelian $p$-groups are considered in Section 2;
we deal only with the case of the prime $p=2$, but the same
method works for odd primes.
The results for the case when
$G$ has order $2$ provide an easy answer
to a question posed in \cite[page 79]{DM}.  
Section 3 contains the analogous theory for actions
of a torus and extends a Borsuk-Ulam for $S^1$ actions
due to Fadell, Husseini and Rabinowitz \cite{FHR, FH} to more general
torus actions.
\section{The Borsuk--Ulam theorem}
We recall a version of the Borsuk--Ulam theorem.
\Prop{\label{BU}
Let $G$ be a finite group.
Suppose that $W$ is a compact, connected, smooth free $G$-manifold
of dimension $n$ and that $V$ is a finite-dimensional real
$G$-module of dimension $k$, with $V^G=0$.
Let $\zeta$ be
the real vector bundle over the orbit space $W/G$
associated with the representation $V$.

Suppose that
the mod $2$ cohomology Euler class $e(\zeta)$
is non-zero. Then, for any continuous $G$-map $f : W\to V$,
the compact subspace
$$
\Zero (f) =\{ x\in W \st f(x)=0\}
$$
has covering dimension greater than or equal to $n-k$. 
}
\begin{proof}
We paste two copies of $W$ together along the boundary
$\partial W$ to form a closed free $G$-manifold
$M = W\cup_{\partial W} W$,
then form the orbit manifold $\bar M= M/G$ and
the associated vector bundle $\xi = M\times_G V$ over $\bar M$.
This bundle is the pullback of the bundle $\zeta$ 
over $W/G$ through the folding map $\pi:  M \to W$.

Now the map $f$ determines a section $s$ of $\xi =M\times_G V$:
$s([x])= [x, f(\pi x)]$ for $x\in M$.
Since $\bar M$ is a closed, connected $n$-manifold and 
$e(\xi)\not=0$, by Poincar\'e duality
there is a class $a\in H^{n-k}(\bar M;\,\Ff_2)$
such that $a\cdot e(\xi)=1\in\Ff_2=H^n(\bar M;\,\Ff_2)$.
By the classical Borsuk--Ulam theorem, as formulated for
example in \cite[Proposition 2.7]{CJ}, the cohomology group
$H^{n-k}(\bar Z;\,\Ff_2)$ of the zero-set
$\bar Z=\Zero (s)\subseteq \bar M$ is non-zero.
Hence the covering dimension of $\bar Z$ is greater than or
equal to $n-k$.
It follows that the inverse image $Z$ of $\bar Z$
under the projection $M\to\bar M$
has covering dimension at least $n-k$,
because $Z \to \bar Z$ is a finite cover.
(See, for example, \cite{J} or \cite[Lemma 2.6]{C}.)

The set $Z$ is a union of two copies of the compact
space $\Zero (f)$ intersecting along $\partial W$. 
So the covering dimension of $\Zero (f)$ is equal to the 
covering dimension
of $Z$ and is greater than or equal to $n-k$.
(Compare \cite[Corollary 2.7]{CS}.)
\end{proof}
There is a more general result for any compact Lie group $G$
and an action that is not necessarily free. We give the formulation
in rational cohomology.
\Prop{\label{BU2}
Let $G$ be a compact Lie group.
Suppose that $W$ is a compact, connected, smooth,
$n$-dimensional $G$-manifold and that $W$ admits an orientation which 
is invariant under the action of $G$.

Let $V$ be a finite-dimensional real
$G$-module of dimension $k$ admitting an orientation
that is fixed by $G$, and let $f : W\to V$ be a continuous $G$-map.

Suppose that either {\rm (i)} the image of
the Borel cohomology Euler class $e(V)\in H^k_G(*;\,\Qq )$
in $H^k_G(W;\,\Qq )$ is non-zero or {\rm (ii)} the
restriction of $f$ to the boundary $\partial W$ is
nowhere zero and the associated relative Euler class
in $H^k_G(W,\partial W;\, \Qq )$ is non-zero.

Then the compact subspace
$$
\Zero (f) =\{ x\in W \st f(x)=0\}
$$
has covering dimension greater than or equal to $n-k$. 
}
The Euler classes in the statement are determined by a choice of
orientation for the vector space $V$.
\begin{proof}
Choose a faithful representation $G\leq \U (\Cc^r)$ and 
write $P=\U (\Cc^r,\Cc^{r+N})$ for the complex Stiefel
manifold of isometric linear maps $\Cc^r\into \Cc^{r+N}$.
With the action of $G\leq \U (\Cc^r)$, $P$ is a closed free
$G$-manifold of dimension $m+\dim G$, say.

Again form the $G$-manifold $M=W\cup_{\partial W}W$ 
as a union of two copies of $W$. 
The connected manifold $\bar M=(P\times M)/G$
of dimension $m+n$ is orientable, and
the map $f$ determines a section $s$
of the vector bundle $\xi =(P\times M\times V)/G$
over $\bar M$ with zero-set $\bar Z =(P\times Z)/G$,
where $Z$ is the union of two copies of $\Zero (f)$. 

If $N$ is sufficiently large (so that $P$ approximates the classifying space $EG$ of $G$), 
then under the hypothesis (i)
$e(\xi )\in H^k(\bar M;\,\Qq )$ is nonzero.
We deduce from the Borsuk-Ulam argument that
$H^{m+n-k}((P\times Z)/G;\,\Qq )$ is non-zero.
Since the Grassmann manifold $P/G$ has dimension $m$,
it follows that $H^{n-k+i}(Z;\,\Qq )$
is non-zero for some $i\geq 0$.
Thus, $Z$ and so also $\Zero (f)$ have covering dimension
$\geq n-k$. 

For (ii), let $\bar W$ be the manifold $(P\times W)/G$
with boundary $\partial\bar W =(P\times\partial W)/G$,
and let $\xi$, now, be the vector bundle $(P\times W\times V)/G$
over $\bar W$. The section $s$ of $\xi$ given by $f$
is non-zero on $\partial \bar W$ and the relative
Euler class 
$e(s;\,\partial \bar W )\in H^k(\bar W,\partial \bar W;\,\Qq )$
is non-zero if $N$ is large.
There is a dual class $a\in H^{m+n-k}(\bar W;\,\Qq )$
such that $a\cdot e(s;\,\partial\bar W)=1\in \Qq =
H^{m+n}(\bar W,\partial\bar W;\, \Qq )$.
The restriction of $a$ to the zero-set $\Zero (s)=(P\times \Zero (f))/G
\subseteq \bar W-\partial\bar W$ is non-zero.
So again the covering dimension of $\Zero (f)$ must be at least $n-k$.
\end{proof}
\section{Elementary abelian $2$-groups}
In this section we consider the case in which $G$ 
is a non-trivial elementary abelian
$2$-group $E$, considered as an $\Ff_2$-vector space of dimension 
$l\geq 1$.

Suppose that $0=E_0\subseteq E_1\subseteq \cdots \subseteq E_l=E$
is a flag in $E$, with $\dim_{\Ff_2}E_i=i$.
The dual vector space $E^*$ parametrizes the $1$-dimensional
real representations of $E$. For a real representation $U$ of
$E$ and $\alpha\in E^*$, we write $U^\alpha$ for the 
$\alpha$-summand,
so that $U=\bigoplus_{\alpha\in E^*} U^\alpha$. 
The annihilator in $E^*$ of a vector subspace $F\subseteq E$
is denoted by $F^\o \subseteq E^*$; it is isomorphic to
$(E/F)^*$. Set $E^i=E_{l-i}^\o$, so that $0=E^0\subseteq \cdots
\subseteq E^i \subseteq \cdots\subseteq E^l=E^*$ is a flag in $E^*$.
We introduce the subspaces
$$
U_i =\bigoplus_{\alpha\in E^i,\, \alpha\notin
E^{i-1}} U^\alpha
\text{\ for $1\leq i\leq l$.}
$$
Thus $U=U^E\oplus \bigoplus_{i=1}^l U_i$, where
$U^E$ ($=U^0$) is the fixed subspace.
This decomposition depends, of course, on the choice of the flag
in $E$.

Let us write 
$$
e(U) =\prod_{\alpha\in E^*} \alpha^{\dim U^\alpha} \in S^*(E^*)
=H^*_E(*;\,\Ff_2)
\text{\quad (and $e(0)=1$)}
$$
in the symmetric algebra of $E^*$ or the Borel cohomology ring of $E$.

We shall need an elementary algebraic result.
\Lem{\label{alg}
Let $\aa$ be an ideal in a commutative ring $A$.
and let $u(T),\, v(T)\in A[T]$ be non-zero polynomials with invertible
leading coefficient.
\par\noindent {\rm (i).}
If $a(T)\in A[T]$ is a polynomial such that
$a(T)u(T)\in\aa [T]$, then $a(T)\in\aa [T]$.
\par\noindent {\rm (ii).}
Suppose that  $\deg (u(T))$ $>\deg (v(T))$. 
Then, if $a\in A$, but $a\notin\aa$,
we have $a\cdot v(T)\notin \aa [T] + (u(T))$. 
}
\begin{proof}
By passing to the quotient $A/\aa$ one can reduce to the case
in which $\aa =0$.
\end{proof}
\Lem{\label{euler}
Suppose that $U$ is a real representation of $E$ such that
$U_i\not=0$ for $i=1,\ldots ,l$.
Then $E$ acts freely on 
$$
\tilde X =\prod_{i=1}^l S(U_i)\subseteq 
U=U^E\oplus \bigoplus_{i=1}^l U_i
$$
and the cohomology ring
of the orbit space $X=\tilde X/E$ is
$$
H^*(X;\,\Ff_2) = S^*(E^*)/(e(U_1),\ldots ,e(U_l)).
$$

Suppose, further, that $V$ is an $E$-module with $V^E=0$, and
let $\xi$ be the vector bundle $\tilde X\times_E V$ over $X$.
Then the $\Ff_2$-Euler class
$e(\xi )$ is non-zero if $\dim U_i > \dim V_i$
for $i=1,\ldots ,l$.
}
\begin{proof}
Notice that $U_i$ and $V_i$  for $i\leq j$ are $E/E_{l-j}$-modules.

Suppose that $v\in E$ and
$x=(x_1,\ldots ,x_l)\in\tilde X$ with $v\cdot x=x$.
Now $E_{l-i+1}/E_{l-i}$ acts freely on $S(V_i)$ for $i=1,\ldots ,l$. So, since $v\cdot x_i=x_i$,
if $v\in E_{l-i+1}$, we can conclude that $v\in E_{l-i}$.
The deduction that $v=0$ is achieved in $l$ steps.
 
The $\Ff_2$-cohomology ring $H^*(X)=H^*_E(\tilde X)$
is calculated step-by-step using the long exact sequences
$$
\cdots\to
H^*_{E/E_{l-j}}(\prod_{i=1}^{j-1} S(U_i)) \Rarr{\cdot e(U_j)}{} 
H^*_{E/E_{l-j}}(\prod_{i=1}^{j-1}S(U_i) ) \to
H^*_{E/E_{l-j}}(\prod_{i=1}^j S(U_i)) \to\cdots
$$
in Borel cohomology.

At the $j$th step, starting with $j=1$,
we apply Lemma \ref{alg} with $A=S^*(E^{l-j+1})$ and
$\aa=(e(U_1), \ldots ,e(U_{j-1}))$.
Choose $\alpha\in E^{l-j}$, $\alpha\notin E^{l-j-1}$,
so that we can identify $S^*(E^{l-j})$ with
the polynomial ring $A[T]$ on $T=\alpha$.
Take $u(T)=e(U_j)$ and $v(T)=e(V_j)$.

Part (i) of Lemma \ref{alg} establishes the injectivity of multiplication by $e(U_j)$ in the exact sequence
and so calculates $H^*_E(\tilde X;\,\Ff_2)=
H^*(X;\,\Ff_2)$.

Part (ii), with $a=e(V_1)\cdots e(V_{j-1})\in A$, so that
$v(T)=a\cdot e(V_j)\in A[T]$, gives the non-vanishing of $e(\xi )$.
\end{proof}
\Ex{(Compare \cite[Theorem 1.2]{G} and \cite[Section 2.7]{BBZ}.)
Suppose that $\phi : \tilde X \to \Rr^m$ is a continuous map and that 
$\dim U_i > m2^{i-1}$ for $i=1,\ldots ,l$.
Then $\phi$ is constant on some $E$-orbit in $\tilde X$.
}
\begin{proof}
Writing $\Rr[E]$ for the group ring of $E$,
set $V= (\Rr [E]/\Rr )^m$. Then $\dim V_i=m2^{i-1}$.
The $E$-map $\tilde X \to V$: 
$x\mapsto [\,\sum_{e\in E} \phi (ex)e\,]$
determines a section of $\xi$. Since
$e(\xi )\not=0$, the section must have a zero,
and this zero in $X$ is the required $E$-orbit.
\end{proof}
As a first application of Lemma \ref{euler} we deduce some results
about group actions on Stiefel manifolds.
\Lem{\label{stiefel}
Let $P$ and $Q$ be finite dimensional Euclidean $E$-modules
with $P^E=0$, $\dim P_i=1$ for $i=1,\ldots ,l$, and $Q^E=0$.
Then, for a given integer $n > l$, 
the group $E$ acts freely on the Stiefel manifold 
$\tilde Y= \O (P,\Rr^n)$
of isometric linear maps $P\into \Rr^n$. 
Let $\eta$ be the
real vector bundle over the orbit space $Y=\tilde Y/E$
associated with the representation $Q$.

If $\dim Q_i \leq n-i$ for each $i=1,\ldots ,l$,
then the Euler class $e(\eta )$ is non-zero.
}
\begin{proof}
Set $U=\Hom(P, \Rr^n)$, so that $U_i=\Hom(P_i,\Rr^n)$ has dimension
$n$.
Let 
$$
h=(h_{i,j}): U \to R\, =\bigoplus_{1\leq i<j\leq l} P_i^*\otimes P_j^*
$$
be given by the inner product on $\Rr^n$: $h_{i,j}(u_1,\ldots ,u_l)
=\langle u_i,u_j\rangle$. The zero-set of $h$ restricted to
$\tilde X$ is exactly the Stiefel manifold $\tilde Y=\O (P,\Rr^n)$.
Notice that $R^E=0$ and $\dim R_i=i-1$.

Take $V=Q\oplus R$.
Then $\xi =\alpha\oplus\beta$, where $\alpha$ and $\beta$ are the
vector bundles over $X$ associated with $Q$ and $R$.
By assumption, $\dim U_i=n>\dim V_i=\dim Q_i +i -1$. 
So the Euler class $e(\xi )=e(\alpha)\cdot e(\beta )$ is non-zero. 

Now the $E$-map $h$ determines a section of $\beta$ with zero-set
equal to $Y$.
By the Borsuk-Ulam theory \cite[Proposition 2.7]{CJ},
the restriction, $e(\eta )$, of $e(\alpha )$ to $Y$ is non-zero.
\end{proof}
\Prop{\label{stiefel_a}
{\rm (Compare \cite[Theorem 5.4]{FH}, \cite[Theorem 1.1]
{frick}.)}
Suppose that $P$ and $Q$ are $E$-modules as in Lemma \ref{stiefel}
such that $\dim Q_i\leq n-i$ for $i=1,\ldots ,l$.
Let $f :\O (P,\Rr^n)\to Q$ be a continuous $E$-equivariant map.
Then the zero-set of $f$ is non-empty and
has covering dimension greater than or equal to $ln-l(l-1)/2-\dim Q$.
}
In particular, we may take $Q_i=\Rr^{n-i}\otimes P_i$.
\begin{proof}
We can apply Proposition \ref{BU} to the manifold $W=\O (P,\Rr^n)$
and the 
\hyphenation{re-presentation}
representation $V=Q$ using Lemma \ref{stiefel}.
\end{proof}
\Rem{\label{flag}
The space $Y=\O (P,\Rr^n)/E$ in Lemma \ref{stiefel}
can be identified with the space of flags 
$0=D^0\subset D^1\subset
\cdots \subset D^l\subseteq \Rr^n$, where $D^j$ is an $\Rr$-subspace 
of dimension $j$, in $\Rr^n$:
send the orbit of $a\in \O (P,\Rr^n)$ to the
flag with $D^j =a(P_1\oplus\cdots\oplus P_j)$.
Let $\delta_i$ be the canonical $(n-i)$-dimensional real vector 
bundle over $B$ with fibre at the flag $(D^j)$ the orthogonal
complement of $D^i$.

It is easy to deduce from the lemma that the product of the Euler 
classes 
$$
e(\delta_1)\cdots e(\delta_l)\in H^{ln-l(l+1)/2}(Y;\,\Ff_2)=\Ff_2
$$
is non-zero.
}
\begin{proof}
By consideration of the projection from the space of flags of
length $n$ to the space of flags of length $l$ we see that
it is enough to deal with the case $l=n$.

So with $l=n$, let $\tau_i$ be the line bundle over
$Y$ associated with the representation $P_i$.
Then $\delta_i$ is the direct sum $\bigoplus_{j=i+1}^n \tau_j$.
So $e(\delta_1)\cdots e(\delta_n)= 
t_1^0\cdots t_i^{i-1}\cdots t_n^{n-1}$, where $t_i=e(\tau_i)$.

Taking $Q_i=\Rr^{n-i}\otimes P_i$ in Lemma \ref{stiefel}, we
deduce that $e(\eta )= t_1^{n-1}\cdots t_i^{n-i}\cdots t_n^0$
is non-zero. 

The result follows from the $\SS_n$-symmetry of
$P =P_1\oplus\cdots \oplus P_n$.
\end{proof}
\Ex{It follows that $l$ sections $s_i$ of $\delta_i$, $i=1,\ldots ,l$,
over the flag manifold $Y$ have a common zero. One can write down
an example with exactly one zero. Choose linearly independent
$v_1,\ldots ,v_l$ in $\Rr^n$ and define the value of $s_i$ at $(D^j)$
to be component of $v_i$ in the orthogonal complement of $D^i$
in the decomposition $\Rr^n =D^i\oplus (D^i)^\perp$.
Then $\bigcap_i \Zero (s_i)$ is precisely the flag $(D^j)$ with
$D^j=\Rr v_1\oplus\cdots\oplus \Rr v_j$.
}
\Rem{\label{ring}
(Compare \cite[Section 3]{BK}.)
Similar methods can be used to describe
the cohomology ring of the flag manifold as 
$$
H^*(\O (P,\Rr^n)/E;\,\Ff_2)=
\Ff_2[t_1,\ldots ,t_l]/(e_1,\ldots ,e_l),
$$
where $e_i=\sum_{r_1+\ldots +r_i=n-i+1}
t_1^{r_1}\cdots t_i^{r_i}$.

More symmetrically, we can replace the classes $e_i$ by
$$
\bar e_i=\sum_{r_1+\ldots +r_l=n-i+1}t_1^{r_1}\cdots t_l^{r_l}
=e_i + \sum_{i<j\leq l} a_{i,j}e_j,
$$
where $a_{i,j}=\sum_{s_j+\ldots +s_l=j-i} t_j^{s_j}\cdots t_l^{s_l}$.

The top-dimensional class is represented by 
$\prod_{i=1}^lt_i^{n-i}$.
}
\begin{proof}
Consider the sphere-bundles
$$
S(P_{i}^*\otimes\zeta_i) \to 
\O (P_1\oplus \cdots \oplus P_{i},\Rr^n)
\to \O (P_1\oplus\cdots \oplus P_{i-1},\Rr^n),
$$
for $1\leq i\leq l$, where $\zeta_i$ is an $E$-vector bundle
of dimension $n-i+1$ with 
$P_1\oplus\cdots\oplus P_{i-1}\oplus\zeta_i=\Rr^n$.
The computation is effected by induction using
the $H^*_E$-Gysin sequences of the sphere-bundles
and Lemma \ref{alg}(i).
The generator $e_i$ corresponds to the Euler class
$e(P_i^*\otimes\zeta_i)=w_{n-i+1}(-(P_1\oplus \cdots \oplus P_i))$.

Write, for $j\geq i$ and $j\leq l$,
$$
e_i[j]=\sum_{r_1+\ldots +r_j=n-i+1}t_1^{r_1}\cdots t_j^{r_j},
$$
so that $e_i=e_i[i]$ and
$e_i[j+1]=e_i[j]+t_{j+1}e_{i+1}[j+1]$.
One easily shows by induction on $k$ that
$$
e_i[i+k]=\sum_{i\leq j\leq i+k}\left(
\sum_{s_j+\ldots +s_{i+k}=j-i} t_j^{s_j}\cdots t_{i+k}^{s_{i+k}} 
\right) e_j.
$$
This then establishes the formula for $\bar e_i$.

The generator of $H^{ln-l(l-1)/2}(Y;\,\Ff_2)=\Ff_2$ was identified as
$e(\eta )$ in Remark \ref{flag} above.
\end{proof}
\Rem{(Compare \cite[Corollary 3.4]{frick}.)
At the cost of losing the $\SS_l$-symmetry, suppose given integers
$1\leq n_1\leq \cdots\leq n_i\leq \cdots \leq n_l\leq n$ such that
$n_i\geq i$.
Consider the submanifold 
$$
\tilde Y'=\{ a\in \O (P,\Rr^n)\st D^j=
a(P_1\oplus\cdots\oplus P_j)\subseteq \Rr^{n_j}\subseteq\Rr^n,\, 
j=1,\ldots ,l\}
$$
of $\tilde Y=\O (P,\Rr^n)$ and the quotient 
$Y'=\tilde Y'/E\subseteq Y$. Let $\eta'$ denote the restriction of $\eta$ to $Y'$.
Then the arguments used in Remark \ref{ring} and Lemma \ref{stiefel}
show that
$$
H^*(Y';\, \Ff_2)=\Ff_2[t_1,\ldots ,t_l]/(e'_1,\ldots ,e'_l),
$$
where $e'_i=\sum_{r_1+\ldots +r_i=n_i-i+1} t_1^{r_1}\cdots t_i^{r_i}$,
and that
the Euler class $e(\eta')$ is non-zero if $\dim Q_i\leq n_i-i$ for
$i=1,\ldots ,l$.
}
Our primary application of Lemma \ref{euler} concerns the
covering dimension of the zero-set of an $E$-map
from $U$ to $V$.
\Prop{\label{main}
Suppose that $U$ and $V$ are $E$-modules as in
Lemma \ref{euler} with $\dim U_i > \dim V_i$ for $i=1,\ldots ,l$.
Let $f : U \to V$ be a continuous $E$-map. Then the
zero-set of $f$ contains a compact free $E$-subspace with
covering dimension greater than or equal to $\dim U -\dim V$.
} 
\begin{proof}
We can apply Proposition \ref{BU} to an equivariant tubular
neighbourhood $W= \tilde X \times D(U^E\oplus\Rr^l)$ of $\tilde X$
in $U$, using Lemma \ref{euler}.
\end{proof}
This implies the following result established in
\cite[Theorem 2.1]{MdMdS}.
\Cor{\label{Mz}
Suppose that $U$ and $V$ are $E$-modules with $V^E=0$ and
$\dim U -\dim V>\dim U^E$,
and suppose that $F\leq E$ is a maximal subgroup such
that $\dim U^F-\dim V^F\geq \dim U -\dim V$.

Let $f : U \to V$ be a continuous $E$-map. Then the
zero-set of $f$ restricted to $U^F$ contains a compact free
$E/F$-subspace with
covering dimension greater than or equal to $\dim U -\dim V$.
}
\begin{proof}
Consider
the action of $E/F$ on $U^{F}$ and $V^{F}$ and the restriction of
$f$ to the fixed-points. We have
$\dim U^F -\dim V^F>\dim (U^F)^E=\dim U^E$. To simplify the
notation, let us assume that $F=0$ and make the abbreviation
$d^\alpha =\dim U^\alpha -\dim V^\alpha$. We show that
there is a flag $(E_i)$ such that $\dim U_i > \dim V_i$
for $i=1,\ldots ,l$. The result will then follow from
Proposition \ref{main}.

Suppose that $E_j$ have been constructed for
$j=l,\ldots ,l-i+1$ for some $i$: $1\leq i <l$. By assumption
$\dim U -\dim V > \dim U^{E_{l-i+1}} -\dim V^{E_{l-i+1}}$,
that is,
$$
\sum_{\alpha\in E^*\, : \,\alpha\notin E^{i-1}} d^\alpha
=
\sum_{E^{i-1}\subseteq E'\subseteq E^*\, :\, \dim E'=i}
{\textstyle \left(\sum_{\alpha\in E' \, : \,\alpha\notin E^{i-1}}
d^\alpha\right)}
\, >\, 0.
$$
Choose a subspace $E'$ such that 
$\sum_{\alpha\in E'\, : \,\alpha\notin E^{i-1}} d^\alpha >0$
and take $E_{l-i}$ to be the annihilator of $E'=E^i$.
\end{proof}
\Cor{\label{cor}
Suppose that $V$ is a Euclidean $E$-module with $V^E=0$ and
that $M$ is a closed $E$-manifold of dimension $n$,
such that the fixed submanifold
$M^E$ is non-empty and has some component of codimension 
strictly greater than $\dim V$. 
Then the zero-set of any $E$-equivariant map $f: M \to V$ contains a compact subspace disjoint from $M^E$ and of covering dimension
at least $n-\dim V$.}
\begin{proof}
Choose a point $x\in M^E$ in a component with codimension
greater than $\dim V$. Take $U$ to be the tangent space $\tau_x M$ 
at $x$ embedded as an 
open $E$-subspace $U\into M$ using the exponential map given by an
$E$-equivariant Riemannian metric on $M$ (mapping $v\in U$
to $\exp_x(\epsilon v/\sqrt{1+\| v|^2})$ for small $\epsilon >0$). 
Now apply Corollary \ref{Mz}.
\end{proof}
We finish this discussion of elementary abelian $2$-groups
with an example involving an infinite dimensional mapping space.
\Prop{\label{loop}
Let $U$ and $V$ be finite-dimensional Euclidean 
$E$-modules such that $U\not=0$,
$\bigcap_{\alpha \, :\, U^\alpha\not=0} \ker \alpha =0$
and $V^E=0$.  

Consider the $E$-space $\map_*(S(\Rr\oplus U),\Rr)$
of real-valued functions on the sphere $S(\Rr\oplus U)$ that are zero at the basepoint $(1,0)$.
Suppose that $f : \map_*(S(\Rr\oplus U),\Rr ) \to V$ is an
$E$-equivariant map. 

Then, for any $d\geq 0$, the zero-set of
$f$ contains a compact free $E$-subspace with covering dimension
greater than or equal to $d$.
}
\begin{proof}
We can choose $\alpha_1,\ldots ,\alpha_l\in E^*$, one at a time,
such that the intersection $\bigcap_{j=1}^i \ker \alpha_j$ has
dimension $l-i$ for $i=1,\ldots ,l$
and then form the flag $(E_i)$ in $E$ such that 
$E_i^\o =\bigcap_{j=1}^i \ker \alpha_j$. With respect to this
flag, which we now fix, each $U_i$ is non-zero.
 
It is easy to see that, for $j\geq 1$, the space 
$P[j]=S^{2j-1}(U^*)$ of
homogeneous polynomial functions $U \to \Rr$ of odd degree $2j-1$
satisfies $\dim P[j]_i \geq \dim U_i$ for $i=1,\ldots ,l$.

For $k\geq 1$, we can embed the $E$-module 
$U[k]=\bigoplus_{j=1}^k P[j]$ in
the mapping space as a space of homogeneous polynomials 
(restricted to $S(\Rr\oplus U)\subseteq \Rr\oplus U$)
$$
\sum_{j=1}^k t^{2k-2j}P[j]\subseteq  \map_*(S(\Rr\oplus U),\Rr )
$$ 
of degree $2k-1$,
where $t$ is the coordinate function on $\Rr$. 

Since $U[k]_i \geq k\dim U_i$, the assertion follows by
applying Proposition \ref{main}, for $k$ sufficiently large,
to the module $U[k]$ instead of $U$.
\end{proof}
\Rem{In the special case that $l=1$, so that $E$ has order $2$,
and $U$ is the non-trivial $1$-dimensional $E$-module $\Rr$ with
the involution $-1$, $\map_*(S(\Rr\oplus U),\Rr )$ is the loop space $\Omega (\Rr, 0)$ with the involution that reverses loops.
Since we can include $\Omega (\Rr ,0)$ in $\Omega (S^n,*)$ 
for any $n\geq 1$,
Proposition \ref{loop} answers a question posed at the end of 
\cite{DM} (although the method 
contradicts the assertion in Proposition 1 of that paper).
}
\Sect{Tori}
Let $L$ be a free abelian group of dimension $l \geq 1$, and
write $T=(\Rr\otimes L)/L$ for the associated torus of rank $l$.
The $1$-dimensional complex representations of $T$ are
parametrized by $L^*=\Hom_\Zz (L,\Zz )$. 
Write $E=\Qq\otimes L$ and $E^*=\Hom_\Qq (E,\Qq )=\Hom_\Zz (L,\Qq )$.
The $\Qq$-projective space on
$E^*$ will be written as $P_\Qq (E^*)$.
For $\lambda\in P_\Qq (E^*)$,
write $L_\lambda =\ker\{ \alpha : L \to \Qq\}$, where $\alpha$
generates $\lambda$. Then $L_\lambda \subseteq L$ has
dimension $l-1$, $T_\lambda = (\Rr\otimes L_\lambda)/L_\lambda$
is a subtorus of $T$ of rank $l-1$ and the quotient
$\Tt_\lambda =T/T_\lambda$ is a circle.

Let $U$ and $V$ be finite-dimensional complex $T$-modules
(which we may assume to be equipped with a $T$-invariant 
Euclidean inner product). 
For $\alpha\in L^*$,
we write $U^\alpha$ and $V^\alpha$ for the $\alpha$-summands
of $U$ and $V$, and, for $\lambda\in P_\Qq (E^*)$,
we set
$$
U_\lambda =\bigoplus_{\alpha\in\lambda-\{ 0\}}U^\alpha=U^{T_\lambda} ,
\quad
V_\lambda =\bigoplus_{\alpha\in\lambda-\{ 0\}} V^\alpha=V^{T_\lambda}\, .
$$
\Prop{Suppose that $U^T=0$ and $V^T=0$.
Then there is a $T$-map $f: U\to V$ with $\Zero (f)=\{ 0\}$ 
if and only if,
for each $\lambda\in P_\Qq (E^*)$ such 
that $U^\lambda\not=0$, there is a $\Tt_\lambda$-map
$S(U_\lambda )\to S(V_\lambda )$.
}
\begin{proof}
(See \cite[Remark 4.4]{BMS0}.)
A $T$-map $S(U)\to S(V)$ restricts on the subspaces
fixed by $T_\lambda$ to a $\Tt_\lambda$-map $S(U_\lambda)\to 
S(V_\lambda)$.
Conversely, given $\Tt_\lambda$-maps $f_\lambda :
S(U_\lambda )\to S(V_\lambda )$ for each $\lambda$ such that
$U_\lambda\not=0$, we can write down a $T$-map
$f: S(U)\to S(V)$:
$$\textstyle
f(\sum_\lambda t_\lambda u_\lambda )=
\sum_\lambda t_\lambda f_\lambda (u_\lambda )
\quad
\text{for $t_\lambda \in [0,1],\, u_\lambda\in S(U_\lambda ),\, 
\sum_\lambda t^2_\lambda =1$}
$$
as the join of the maps $f_\lambda$.
\end{proof}
\Ex{Take $l=1$ and $L=\Zz$, so that $T$ is the circle $\Rr /\Zz$.
Let $a,\,b,\, c\geq 1$ be positive integers with $a$ and $b$ coprime.
Take $U=\Cc\oplus\Cc$ with $t+\Zz$ acting as $(x,y)\mapsto
(\e^{2\pi\i act}x,\e^{2\pi\i bct}y)$ and $V=\Cc\oplus\Cc$ with
the action $(x,y)\mapsto (\e^{2\pi\i abc t}x,\e^{2\pi\i ct}y)$.
Choose $a',\, b'\geq 1$ with $aa'-bb'=1$.
Then $f: U\to V$
$$
(x,y)\mapsto (x^b+y^a, x^{a'}\bar y^{b'})
$$
is $T$-equivariant and $\Zero (f)=\{ 0\}$.
}
We write
$$
e(U) =\prod_{\alpha\in L^*} \alpha^{\dim U^\alpha} \in 
S^*(E^*)=H^*_T(*;\,\Qq )
$$
for the Euler class in $T$-equivariant Borel cohomology.

\Prop{{\rm (Compare \cite[Theorem (2.3)]{FHR}.)}
Suppose that $U^T=0$, $V^T=0$, $\dim U >\dim V$ and $\phi : U\to\Rr$
is a continuous $T$-map such that $\phi (0)<0$ and $\phi (x)>0$
for $\| x\|$ sufficiently large.

If $f : U\to V$ is a $T$-map,
then the intersection $\Zero (f)\cap\Zero (\phi )$
has covering dimension greater than or equal
to $2(\dim_\Cc U-\dim_\Cc V)-1$.
}
\begin{proof}
(Compare \cite[Section 5]{BMS}.)
Choose $v\in L$ such that $\alpha (v)\not=0$ for all
$\alpha\in L^*-\{ 0\}$ such that $U^\alpha\not=0$
or $V^\alpha\not=0$. Let $\rho : \Tt =\Rr /\Zz \to
T=(\Rr\otimes L)/L$ be the associated homomorphism
$t+\Zz\mapsto tv+L$. Then $U^{\rho (\Tt )}=0$
and $V^{\rho (\Tt )}=0$.

Choose radii $r$ and $R$, $0<r<R$, such that $\phi (x)<0$ if
$\| x\| =r$ and $\phi (x)>0$ if $\| x\| =R$.
The annulus $W=\{ x\in U\st r\leq \| x\| \leq R\}$ is a compact
manifold of dimension $n=2\dim_\Cc U$ which is
$T$-equivariantly diffeomorphic to $D(\Rr )\times S(U)$.

Write $k=2\dim_\Cc V+1$. Now apply Proposition \ref{BU2} with
condition (ii) and $G=\Tt$ to the map $(\phi ,f): W \to \Rr\oplus V$.
The relative Euler class in $H^k_\Tt (W,\partial W;\,\Qq )
=H^{k-1}_\Tt (S(U);\,\Qq )=\Qq$ is the image of $e(V)$ and is non-zero
(because $k-1=2\dim_\Cc V\leq 2(\dim_\Cc U-1)$).
\end{proof}
Suppose that $0=E_0\subseteq \cdots\subseteq E_i \subseteq\cdots
\subseteq E_l=E$, with $\dim_\Qq E_i=i$, is a flag in the $\Qq$-vector space $E$.
Put $L_i=L\cap E_i$ and $T_i=(\Rr\otimes L_i)/L_i\leq T$.
Writing $E^i=E^\o_{l-i}\subseteq E^*$ for the annihilator
of $E_{l-i}$, we define
$$
U_i=\bigoplus_{\lambda \in P_\Qq (E^i),\, \lambda\notin
P_\Qq (E^{i-1})} U_\lambda.
$$
So $U=U^T\oplus\bigoplus_{i=1}^l U_i$.

\Lem{\label{eulerT}
Suppose that $U$ is a complex representation of $T$
such that $U_i\not=0$ for all $i$. Then $T$ acts with finite isotropy
groups on
$$
\tilde X=\prod_{i=1}^l S(U_i) \subseteq U=U^T\oplus
\bigoplus_{i=1}^l U_i
$$
and $H^*_T(\tilde X;\, \Qq )=S^*(E^*)/(e(U_1),\ldots ,e(U_l))$.

If $V$ is a complex $G$-module with $V^T=0$ and $\dim U_i >\dim V_i$ 
for $i=1,\ldots ,l$, 
then the image of $e(V)$ in $H^*_T(\tilde X;\,\Qq )$ is non-zero.
}
\begin{proof}
Consider $x=(x_1,\ldots ,x_l)\in \tilde X$ fixed by $v+L\in T$,
where $v\in \Rr\otimes L$. 
For each $i$, there is some
$\alpha_i\in L^*$ with $\alpha_i\in E^i=E_{l-i}^\o$,
$\alpha_i\notin E^{i-1}=E_{l-i+1}^\o$ such that the component
of $x_i$ in $U^{\alpha_i}$ is non-zero.
Since $v+L$ fixes $x_i$, we have $\alpha_i(v)\in\Zz$.
Because the $\alpha_i$ form a $\Qq$-basis of $E^*$, 
the vector $v\in \Rr\otimes L$ lies in $E=\Qq\otimes L$
and $v+L\in T$ has finite order.
Thus the isotropy group of $x$ is finite.

The rest of the proof follows
{\it mutatis mutandis} the argument in Lemma \ref{euler}
using long exact sequences
$$
\cdots\to
H^*_{T/T_{l-j}}(\prod_{i=1}^{j-1} S(U_i)) \Rarr{\cdot e(U_j)}{} 
H^*_{T/T_{l-j}}(\prod_{i=1}^{j-1}S(U_i) ) \to
H^*_{T/T_{l-j}}(\prod_{i=1}^j S(U_i)) \to\cdots
$$
in Borel cohomology. 
\end{proof}
There are complex analogues of Propositions \ref{stiefel_a}
and \ref{main}.
\Prop{Let $P$ and $Q$ be finite dimensional Hermitian $T$-modules
with $P^T=0$, $Q^T=0$, and $\dim_\Cc P_i=1$ for $i=1,\ldots ,l$.
Suppose that, for some $n > l$,
$f: \U (P,\Cc^n)\to Q$ is a $T$-equivariant
map from the complex Stiefel manifold of isometric $\Cc$-linear maps
$P\into \Cc^n$ to $Q$ and that $\dim_\Cc Q_i \leq n-i$ for each 
$i=1,\ldots ,l$.
Then the zero-set of $f$ is a non-empty free $T$-space with covering
dimension at least $2ln-l^2-2\dim_\Cc Q$.
}
\begin{proof}
This can be established by using
Proposition \ref{BU2} with condition (i) taking $G=T$,
$W=\U (P,\Cc^n)$ and $V=Q$.
\end{proof}
\Prop{\label{mainT}
Suppose that $U$ and $V$ are $T$-modules as in
Lemma \ref{eulerT} with $\dim U_i > \dim V_i$ for $i=1,\ldots ,l$.
Let $f : U \to V$ be a continuous $T$-map. Then the
zero-set of $f$ contains a compact $T$-subspace with
finite isotropy groups and
covering dimension greater than or equal to 
$2(\dim_\Cc U -\dim_\Cc V)$.
} 
\begin{proof}
We can apply Proposition \ref{BU2} with condition (i) 
to an equivariant tubular
neighbourhood $W= \tilde X \times D(U^E\oplus\Rr^l)$ of $\tilde X$
in $U$, using Lemma \ref{eulerT}, with $k=2\dim_\Cc V$ and
$n=2\dim_\Cc U$. The orientations, being determined by complex structures, are invariant under the action of $G=T$.
\end{proof}
\par\noindent
{\bf Acknowledgment}
The idea of estimating the dimension of the free part of the zero-set
can be found in \cite{DM}, and I am grateful to Professor Miklaszewski
for correspondence about his work.

\end{document}